# Lattices and the Geometry of Numbers

Sourangshu Ghosh[a],

[a]Undergraduate Student, Department of Civil Engineering,
Indian Institute of Technology Kharagpur, West Bengal, India

## 1. ABSTRACT

In this paper we discuss about properties of lattices and its application in theoretical and algorithmic number theory. This result of Minkowski regarding the lattices initiated the subject of Geometry of Numbers, which uses geometry to study the properties of algebraic numbers. It has application on various other fields of mathematics especially the study of Diophantine equations, analysis of functional analysis etc. This paper will review all the major developments that have occurred in the field of geometry of numbers. In this paper we shall first give a broad overview of the concept of lattice and then discuss about the geometrical properties it has and its applications.

## 2. LATTICE

Before introducing Minkowski's theorem we shall first discuss what is a lattice.

***Definition 1:*** A lattice $\tau$ is a subgroup of $\mathbf{R^n}$ such that it can be represented as

$$\tau = a_1 \mathbf{Z} + a_2 \mathbf{Z} + \ldots + a_m \mathbf{Z}$$

Here $\{a_i\}$ are linearly independent vectors of the space $\mathbf{R^n}$ and $m \leq n$. Here $\mathbf{Z}$ is the set of whole numbers.

We call these vectors $\{a_i\}$ the *basis* of the lattice. By the definition we can see that a lattice is a subgroup and a free abelian group of rank m, of the vector space $\mathbf{R^n}$. The rank and dimension of the lattice is $m$ and $n$ respectively, the lattice will be complete if $m = n$. This definition is not only limited to the vector space $\mathbf{R^n}$. It can be extended to any arbitrary field $\mathbf{F}$, in which the basis vectors $\{a_i\}$ will belong to the field $\mathbf{F}$. In this article we shall discuss about both complete or full-rank lattices and incomplete lattices. We will now define another terms known as the fundamental mesh.

***Definition 2:*** The set of elements which can be denoted as

$$\varphi(a) = \{a_1 \boldsymbol{v_1} + a_2 \boldsymbol{v_2} + \ldots + a_m \boldsymbol{v_m} | v_i \in R, 0 \leq v_i \leq 1\}$$

is called a fundamental mesh of the lattice.

A very important thing to notice is that not every given set of vectors $\{a_i\}$ forms the basis of a given lattice $\tau$. In the next lemma we shall state the condition for a given set of vectors $\{a_i\}$ to form the basis of a given lattice $\tau$.

***Lemma 1:*** A given set of vectors $\{a_i\}$ form the basis of a given lattice $\tau$ if it satisfies the following condition:

$$\varphi(a) \cap \tau = \{\mathbf{0}\}$$

***Proof:*** *The* lattice $\tau$ is the set of all their integer combinations of the basis vectors $\{a_i\}$. We also know that $\varphi(a)$ is the set of linear combinations of basis vectors $\{a_i\}$ with the coefficients $v_i \in R, 0 \leq v_i \leq 1$. Therefore the only element in common is the zero vector $\mathbf{0}$.

Let us now state a lemma regarding the properties of the lattice $\tau$. This proof can be found in many articles and books. We state here the elegant proof as stated and given by Comeaux[1] and Neukirch[2]

# Lattices and the Geometry of Numbers

***Lemma 2:*** A lattice $\tau$ in $V$ is complete if and only if there exists a bounded subset $M \subseteq V$ such that the collection of all translates $M + \gamma, \gamma \in \tau$ covers the whole space $V$.

***Proof:*** Comeaux[1] proves the theorem by proving it in the other direction. Let us first assume that $M$ is a bounded subset such that the collection of all translates $M + \gamma, \gamma \in \tau$ covers the whole space $V$. Let us also denote the subspace spanned by $\tau$ as $P$. Now to prove this theorem, we have to prove $P = V$ or in other words if we can show that every element $v \in V$, also belongs to $P$ we are done. Now we note that $V = \cup_{\gamma \in \tau}(M + \gamma)$, we also note that $v \in V$ and $a_v \in M$ and $\gamma_v \in \tau$. But $\tau$ is itself a subset of $P$. We can therefore write the following expression:

$$vv = a_v + \gamma_v$$

Dividing both sides of the equation by $v$ and taking the limit of $v$ as infinity, we get

$$v = \lim_{v \to \infty}(\frac{a_v}{v}) + \lim_{v \to \infty}(\frac{\gamma_v}{v}) = \lim_{v \to \infty}(\frac{\gamma_v}{v}) \in V'$$

We will now define another term known as the determinant of the lattice $\tau$, or the volume of the parallelepiped spanned by the fundamental mesh. Before going to the technical definition let us discuss a more intuitive way to understand the notion of determinant of the lattice $\tau$, the determinant of the lattice $\tau$ can be expressed as:

$$\det(\tau) = \lim_{k \to \infty} \frac{vol(k)}{Number\ of\ lattice\ points\ inside}$$

Where $k$ is the radius of the ball.

***Definition 3:*** The determinant of the lattice $\tau$, or the volume of the parallelepiped spanned by the fundamental mesh spanned by the basis vectors $\{a_i\}$ is given as

$$\det(\tau) = \sqrt{\det(A^T A)}$$

If the lattice if full-rank, the expression reduces to

$$\det(\tau) = |\det(A)|$$

We have discussed above what a lattice is. An important thing to notice is that it is not necessary that two distinct basic vectors $\{a_i\}$ and $\{a_i'\}$ shall span distinct lattices. The span of two different basic vectors $\{a_i\}$ and $\{a_i'\}$ can give us the same lattice. We shall now state another lemma regarding this question as given by Oded Regev[3].

***Lemma 3:*** The span of two different basic vectors $\{a_i\}$ and $\{a_i'\}$ will give us the same lattice if and only if $A' = AU$, where $U$ is a uni-modular matrix.

***Proof:*** A uni-modular matrix is matrix whose determinant is $\pm 1$.. We first assume that the span of two different basic vectors $\{a_i\}$ and $\{a_i'\}$ is equal. From this point onwards we shall use $\{a_i\}$ and $A$ interchangeably. Now note that each column of the matrix $A$ which is $a_i$ belongs to the span of $A'$, which tells us that there must exist a uni-modular U for which $A' = AU$. Similarly we can otherwise say that there must exist a uni-modular $U'$ for which $A = A'U'$. Substituting the later expression into the first one we get

$$A' = AU = A'U'U$$

Therefore we can also say that

$$A'^T A' = (U'U)'(A'^T A')U'U$$

Now if we take determinant on both sides of the equation we get

$$det(A'^T A') = det(U'U)^2 det(A'^T A')$$

This gives us $det(U'U)^2 = 1$ i.e $det(U'U) = \pm 1$ or in other words $det(U')det(U) = \pm 1$. Now notice that both $U'$ and $U$ is integer matrices and hence the determinant shall also be integers. Therefore we get

$$det(U) = \pm 1$$

We shall now discuss about the method of Gram-Schmidt Orthogonalization. This again can be found in many textbooks but for the present purposes we shall be stating it as done by Oded Regev[3]. It is a method of ortho-normalizing (such that the new set of basis vectors $\tilde{a}_i$ satisfy $\langle \tilde{a}_i, \tilde{a}_j \rangle = 0$) a given set of linearly independent basis vectors $\{a_i\}$, which may be a part of an inner product space equipped with a given inner product.

***Definition 4:*** For a sequence of $n$ linearly independent basis vectors $\{a_i\}$. We now define the Gram-Schmidt orthogonalization as the sequence of vectors $\{\tilde{a}_i\}$ defined by

$$\tilde{a}_i = a_i - \sum_{j=1}^{i-1} \mu_{i,j} \tilde{a}_j \text{ where } \mu_{i,j} = \frac{\langle a_i, \tilde{a}_j \rangle}{\langle a_j, \tilde{a}_j \rangle}$$

One should note that $\{\tilde{a}_i\}$ do not necessarily form a basis of $\tau$ and the order of the linearly independent basis vectors $\{a_i\}$ matters due to which we should consider it as a "*sequence of vectors*" rather than a "*set of vectors*". Geometrically speaking in this method the basis vector $\tilde{a}_i$ is defined as difference between the basis vector $a_i$ and its projection onto the subspace which is generated by $a_1, a_2, \ldots, a_{i-1}$ which shall be by definition same as the subspace generated by $\tilde{a}_1, \tilde{a}_2, \ldots, \tilde{a}_{i-1}$. By subtracting the projection from the vector we are making sure that it will be orthogonal to subspace generated by $\tilde{a}_1, \tilde{a}_2, \ldots, \tilde{a}_{i-1}$.

Let us now state another important property of the Gram-Schmidt Orthogonalization due to Oded Regev[3]. Let there be n linearlly independent vectors denoted as $\{a_i\}$. Therefore we get the ortho-normal basis vectors given by Gram-Schmidt Orthogonalization as $\{\tilde{a}_i/\|\tilde{a}_i\|\}$. Then we can represent the new basis vectors $\{\tilde{a}_i\}$ in terms of basis vectors $\{a_i\}$ as given as the columns of the $m \times n$ matrix.

$$\begin{bmatrix} \|\tilde{a}_1\| & \mu_{2,1}\|\tilde{a}_1\| & \ldots & \mu_{n,1}\|\tilde{a}_1\| \\ 0 & \|\tilde{a}_2\| & & \mu_{n,2}\|\tilde{a}_2\| \\ \vdots & & & \\ 0 & \ldots & 0 & \|\tilde{a}_n\| \\ 0 & \ldots & 0 & 0 \\ \vdots & & \vdots & \vdots \\ 0 & \ldots & 0 & 0 \end{bmatrix}$$

This matrix will be a full-rank lattice if $m = n$, in which case the matrix will become a upper-triangular matrix. The volume of the parallelopoid $P(\{a_i\})$ will be equal to $det(L(\{a_i\}))$ which in turn is same as $\prod_{i=1}^{n} \|\tilde{a}_i\|$.

### 3. SUCCESIVE MINIMA

Having discussed about the lattice and its properties, we shall now discuss about the concept of Minkowski Theory. Before that we shall state some definitions regarding the properties of subset.

***Definition 5:*** A subset $X \in V$ is defined to be as ***centrally symmetric***, if for every point $x \in X$, $-x \in X$ also holds

***Definition 6:*** A subset $X \in V$ is defined to be as ***convex***, if for every two point $x, y \in X$, the line segment defined as

# Lattices and the Geometry of Numbers

$$\{ty + (1-t)x | 0 \leq t \leq 1\}$$

Also is contained inside the set $V$.

One very important parameter of interest in the study of lattices is the determination of the length of the shortest nonzero vector in the lattice which we call as $\lambda_1$. The other successive minima distances are written as, $\lambda_2,..,\lambda_n$. We shall now study more about the concept of successive minima of a symmetric convex body and the successive inner and outer radii of the body as stated by Henk[4,5] and Oded Regev[3]. But first we shall give a definition of $i^{th}$ successive minima of a lattice $\tau$ of rank n.

***Definition 7:*** We define the $i^{th}$ successive minima of a lattice $\tau$ of rank n as

$$\lambda_i(\tau) = \inf\{r | \dim(span(\tau \cap \overline{B}(0,r))) \geq i\}$$

Here $\overline{B}(0,r)$ is a closed ball of radius $r$ around the point 0.

By the definition itself it is clear that we shall have $\lambda_1 \leq \lambda_2 \leq .. \leq \lambda_n$, also we here note that the first successive minima shall be greater than the minimum of all the sequence of orthogonal basis vectors that span $\tau$ that we found earlier by Gram-Schmidt orthogonalization i.e. $\lambda_1 \geq \min_{i=1,2..,n} \tilde{a}_i$. We shall state and prove this observation in our next theorem due to by Micciancio[6].

***THEOREM 1:*** Let $A$ be the basis of the full-rank lattice $\tau$ of rank n, and let $\tilde{A}$ be its Gram-Schmidt orthogonalization. Then the following inequality holds:

$$\lambda_1(\tau) \geq \min_{i=1,2..,n} ||\tilde{a}_i|| > 0$$

***Proof:*** Let us consider the lattice vector $Ax$, and let $k$ be the biggest index such that $x_k$ is not zero. Now notice that to prove the theorem it is sufficient to prove

$$||Ax|| \geq ||\tilde{a}_k|| \geq \min_{i=1,2..,n} ||\tilde{a}_i||$$

Now we take the scalar product of our lattice vector and $\tilde{a}_i$. Now because of the orthogonality of $\tilde{a}_k$ and $a_i$ (for $i < k$) we can write

$$\langle Ax, \tilde{a}_k \rangle = \sum_{i \leq k} \langle a_i x_i, \tilde{a}_k \rangle = x_k \langle a_k, \tilde{a}_k \rangle = x_k ||\tilde{a}_k||^2$$

Now by Cauchy-Schwarz Inequality we know that $||Ax||.||\tilde{a}_k|| \geq |\langle Ax, \tilde{a}_k \rangle| \geq x_k ||\tilde{a}_k||^2$. Using $|x_k| \geq 1$ we get

$$||Ax|| \geq ||\tilde{a}_k||$$

We shall now state the famous Minkowski's Lattice Point Theorem or Minkowski's Convex Body Theorem as stated in by Comeaux[1], Neukirch[2] and Matousek[7] and Minkowski[8]. Before stating Minkowski's theorem we shall go in the traditional way to prove the theorem, that is to start with a theorem of Blichfeld, and try to prove Minkowski's theorem as a corollary of a theorem of Blichfeld as done by Micciancio[6]. These will later help us derive upper bounds on the successive minima. We shall give a concise proof for this 2 theorems, for a more detailed proof please see the above mentioned references.

***Theorem 2:*** (**Blichfeld Theorem**) For a given lattice $\tau$ and a subset $X$ which belong to span($A$), if $Vol(X) > \det(\tau)$ there must exist two distinct points let's say $z_1$ and $z_2$ which belong to $X$. Such that $z_1 - z_2$ also belong to $\tau$.

# Lattices and the Geometry of Numbers

**Proof:** Micciancio[5] proved the theorem by first considering the collection of sets defined as follows

$$X_y = X \cap (y + P(A))$$

Here $y \in \tau$. By its construction these sets are pair-wise disjoint form a partition of $X$. Therefore we can also say that

$$vol(X) = \sum_{y \in \tau} vol(X_y)$$

Micciancio[5] then proved that the sets $X_y$ can't be all mutually disjoint. Note that $X_y - y$ is a subset of span($A$) and also the fact that $vol(X_y) = vol(X_y - y)$ because the sets $X_y - y = (X - y) \cap \text{span}(A)$ are all subsets of span($A$), we therefore have the following inequality

$$vol(\text{span}(A)) = \det(\tau) < vol(X) = \sum_{y \in \tau} vol(X_y) = \sum_{y \in \tau} vol(X_y - y)$$

The previous inequality means

$$\sum_{y \in \tau} vol(X_y - y) > vol(\text{span}(A))$$

This means that these sets $X_y$ cannot be disjoint as we also have $X_y - y$ as a subset of span($A$). Therefore theres must exist two distinct non-trivial points belonging to the lattice $\tau$ let's say $x, y$ such that the intersection of $(X_y - y)$ and $(X_x - x)$ is not null and therefore contain a point let's say $z$ in it. As it is part of both the sets $(X_y - y)$ and $(X_x - x)$, we can hence write two distinct non-trivial points $z_1$ and $z_2$ such that $z_1 = z + x \in X_x$ and $z_2 = z + y \in X_y$. And as both are subsets of $X$. These two points satisfy

$$z_1 - z_2 = x - y \in \tau$$

We can now prove the Minkowski's Convex Body Theorem as a corollary to Blichfeld Theorem. These will later help us to give an some bounds to the length of the shortest vector possible in a given lattice $\tau$

**Theorem 3:** (**Minkowski's Convex Body Theorem**) Let $\tau$ be a complete lattice in a euclidean vector space $V$ and let $X$ be a centrally symmetric, convex subset of $V$. Suppose the following condition holds:

$$Vol(X) > 2^n \det(\tau)$$

If the above mentioned condition holds then $X$ must contain a nonzero lattice point $\gamma \in \tau$.

**Proof:** Let us take the set $\hat{X} = \frac{1}{2}X = \{x | 2x \in X\}$. Note that the volume of $\frac{1}{2}X$ satisfies $vol\left(\frac{1}{2}X\right) = 2^{-n} vol(X) > \det \tau$. By earlier theorem 1 by Blichfeld there must exist two distinct points let's say $z_1$ and $z_2$ which belong to $\frac{1}{2}X$. Such that $z_1 - z_2$ also belong to $\tau$ excluding the trivial point 0. Now as $z_1$ and $z_2$ which belong to $\frac{1}{2}X$. $2z_1$ And $2z_2$ shall belong to $X$. Now we assumed that $X$ is convex and symmetric, therefore $-2z_2$ shall belong to $X$.

$$z_1 - z_2 = \frac{2z_1 - 2z_2}{2} \in X$$

This point is a non-trivial lattice point contained in the set $X$.

We shall now give some corollaries to the Minkowski's Convex Body Theorem. The first one is due to Oded Regev[3].

# Lattices and the Geometry of Numbers

***Corollary 1:*** The volume of an n-dimensional ball of radius r is $vol(\boldsymbol{B}(\boldsymbol{0}, \boldsymbol{r})) \geq (\frac{2r}{\sqrt{n}})^n$

***Proof:*** The n-dimensional ball of radius $\boldsymbol{r}$ must contain inside completely the cube whose side length shall be twice of $\frac{r}{\sqrt{n}}$. Or in other words the volume of the n-dimensional ball of radius $\boldsymbol{r}$ must be at least $(\frac{2r}{\sqrt{n}})^n$.

***Corollary 2:*** **(Minkowski's First Theorem)** For any full-rank lattice $\boldsymbol{\tau}$ of rank $n$, the following inequality holds:

$$\lambda_1 \leq \sqrt{n}(det(\boldsymbol{\tau}))^{1/n}$$

***Proof:*** This can be easily proved once we notice that the ball $\boldsymbol{B}(\boldsymbol{0}, \lambda_1)$ contains no lattice point, therefore we can write according to the Minkowski convex body theorem we defined earlier the following:

$$vol(\boldsymbol{B}(\boldsymbol{0}, \lambda_1)) \leq 2^n det(\boldsymbol{\tau})$$

Also by corollary 1 we know that

$$vol(\boldsymbol{B}(\boldsymbol{0}, \lambda_1)) \geq (\frac{2\lambda_1}{\sqrt{n}})^n$$

Comparing the two inequalities we get

$$\lambda_1 \leq \sqrt{n}(det(\boldsymbol{\tau}))^{1/n}$$

A important thing to note is that in the discussion above we used the $l_2$ norm only. But we can similarly extend extend Minkowski's theorem to other norms. For that all we have to do is to compute the volume in the given norm of the ball. Let us consider the $l_\infty$ norm for the current discussion. The $l_\infty$ norm of a vector is denoted by $\|x\|_\infty = \max_i |x_i|$. We shall now state a corollary to the Minkowski's theorem as stated by Micciancio[6].

***Corollary 3:*** It can show that for every given lattice $\boldsymbol{\tau}$ of rank $n$ must contain a nonzero vector x such that the following inequality holds:

$$\|x\|_\infty \leq (\det \boldsymbol{\tau})^{1/n}$$

***Proof:*** Let us define the distance $l$ as

$$l = \min\{\|x\|_\infty : x \in \boldsymbol{\tau} \setminus \{\boldsymbol{0}\}\}$$

We shall now assume that $l > (\det \boldsymbol{\tau})^{1/n}$. As earlier done in corollary 2 we shall take the hypercube $\boldsymbol{C} = \{x : \|x\| < l\}$. The volume of the hypercube will be

$$vol(C) = (2l)^n > 2^n \det(\boldsymbol{\tau})$$

Now as $\boldsymbol{C}$ is convex, symmetric we can say using the Minkowski's Convex Body Theorem that there exists a non-trivial lattice point $\boldsymbol{x}$ inside the hypercube. But note that we have $\|x\|_\infty \leq l$ by definition of hypercube $\boldsymbol{C}$ itself. Therefore we arrive at a contradiction to the minimality of $l$.

***Corollary 4[Nguyen[9]]:*** Any $d-$dimensional lattice $\boldsymbol{\tau}$ of $R^n$ contains a nonzero $x$ such that the following inequality holds:

$$\|x\| \leq 2(\frac{vol(\tau)}{v_d})^{1/d}$$

Here $v_d$ denotes the volume of the closed unitary hyperball of $R^d$.

# Lattices and the Geometry of Numbers

In the previous theorem we have seen a bound of shortest nonzero vector which is the first successive minimum $\lambda_1$ by the Minkowski's First Theorem. One important thing to notice is that we could have directly proven corollary 2 using corollary 3 as for any $n$-dimensional vector $x$ the inequality $\|x\| < \sqrt{n}\|x\|_\infty$ holds. We shall now state and prove the Minkowski's second theorem which gives us a bound considering the geometric mean of all the successive minimum of a given lattice $\lambda_i$ as stated by various mathematicians[10-17]. As it is the geometric mean of $\lambda_i$ it is clearly greater than $\lambda_1$, hence this bound can be considered stronger than the first one.

***THEOREM 4:*** **(Minkowski's Second Theorem):** For any full-rank lattice $\tau$ of a given rank $n$ the following inequality holds:

$$\frac{\sqrt{n}(det(\tau))^{1/n}}{n!^{1/n}} \leq (\prod_{i=1}^{n}\lambda_i)^{1/n} \leq \sqrt{n}(det(\tau))^{1/n}$$

***Proof:*** The following short but elegant proof we shall discuss is given by Evertse[18].

Let us take a full-rank lattice $\tau$ of a given rank $n$ having $n$ linearly independent vectors $x_1, x_2, \ldots x_n$. Note that we have earlier proved that $\|x\|_i = \lambda_i$. We shall first prove the left side of the inequality which states that

$$\frac{\sqrt{n}(det(\tau))^{1/n}}{n!^{1/n}} \leq (\prod_{i=1}^{n}\lambda_i)^{1/n}$$

Let us take the convex symmetric body defined by

$$C := \{\sum_{i=1}^{n} x_i \lambda_i^{-1} v_i :, \sum_{i=1}^{n}|x_i| \leq 1\}$$

Now notice that we can get the convex symmetric body $C$ from another convex symmetric body $L$ defined as

$$L := \{x \in R^n : \sum_{i=1}^{n}|x_i| \leq 1\}$$

By the linear transformation $\varphi$ which takes the basis vectors of $R^n$ and gives $\lambda_i^{-1} v_i$. The Jacobian of the transformation will be $det\varphi$. Therefore the volume of the convex symmetric body will be

$$vol(C) = vol(L)det\varphi = (\frac{2^n}{n!}) * \frac{|det(\tau)|}{\prod_{i=1}^{n}\lambda_i}$$

Note that $C$ is the smallest convex set symmetric about the point **0** containing the points $\lambda_i^{-1} v_i$. Therefore the cube whose side length is twice of $\frac{1}{\sqrt{n}}$ must contain the symmetric convex set $C$. Therefore we can write

$$vol\ of\ cube = (\frac{2}{\sqrt{n}})^n \geq vol(C) = (\frac{2^n}{n!}) * \frac{|det(\tau)|}{\prod_{i=1}^{n}\lambda_i}$$

Which can be otherwise be written as $\frac{\sqrt{n}(det(\tau))^{1/n}}{n!^{1/n}} \leq (\prod_{i=1}^{n}\lambda_i)^{1/n}$. We now similarly prove the right hand side of the inequality by taking the convex symmetric body defined by

$$C := \{\sum_{i=1}^{n} x_i \lambda_i^{-1} v_i : x \in R^n : |x_1| \leq 1, \ldots, |x_n| \leq 1\}$$

Now notice that we can get the convex symmetric body $C$ from another convex symmetric body $L$ defined as

$$L := \{x \in R^n : |x_1| \leq 1, \ldots, |x_n| \leq 1\}$$

As we have earlier done we take the linear transformation $\varphi$ which takes the basis vectors of $R^n$ and gives $\lambda_i^{-1} v_i$. The Jacobian of the transformation will be $det\varphi$. Therefore the volume of the convex symmetric body will be

$$vol(C) = vol(L) det\varphi = 2^n * \frac{|det(\tau)|}{\prod_{i=1}^n \lambda_i}$$

Evertse[18] showed that the successive minima of $C$ are same as lattice $\tau$. For a given $\lambda > 0$ we have $\lambda C$ defined as

$$\lambda C := \{\sum_{i=1}^n x_i \lambda_i^{-1} v_i : x \in R^n : |x_1| \leq \lambda, \ldots, |x_n| \leq \lambda\}$$

This formulation implies that if $\lambda = \lambda_i$, Then the defined volume shall contains $i$ linearly independent basis vectors $v_1, \ldots, v_i$. Now let us assume $\lambda$ is lesser than $\lambda_i$, we also take a vector $p = \sum_{j=1}^n p_j v_j$. Now by the formulation itself of the defined volume $\lambda_i C$ we have $|p_j| \leq 1$ for $j = i, \ldots, n$ which implies that $p_j = 0$ for $j = i, \ldots, n$. Therefore all the lattice points shall lie in the volume spanned by the $i - 1$ linearly independent basis vectors $v_1, \ldots, v_{i-1}$. Therefore we can say that $\lambda_i$ is the $i$th successive minima as this space cannot contain $i$ linearly independent points Therefore the cube whose side length is twice of $\frac{1}{\sqrt{n}}$ must be itself contained inside the symmetric convex set $C$ which was the containing volume in the lower limit derivation we have done earlier. Therefore we can write

$$vol\ of\ cube = (\frac{2}{\sqrt{n}})^n \leq vol(C) = 2^n * \frac{|det(\tau)|}{\prod_{i=1}^n \lambda_i}$$

Which can be otherwise be written as $(\prod_{i=1}^n \lambda_i)^{1/n} \leq \sqrt{n}(det(\tau))^{1/n}$.

We shall now define the Hermite constant, which is a fundamental quantity in mathematics that determines how much a lattice element can be in a lattice belonging to Euclidean space.

## 4. SPHERE PACKINGS

A problem of huge interest is to know what fraction of $R^n$ can be covered by equal balls that do not intersect except along their boundaries as stated by Conway and Sloane[19] written in Nguyen[9]. This is an open problem for all values of n greater than 3. This problem can be restated as what is the densest packing that is derived from lattices, which is famously called the lattice packing problem. To address the problem we shall first define the density of the lattice and the hermite constant both of which are important parameters to measure the packing in the lattice.

***Definition 8 [Nguyen[9]]:*** The density $\delta(\tau)$ of the lattice packing is equal to the ratio between the volume of the n-dimensional ball of diameter $\lambda_1(\tau)$ and the volume of any fundamental domain of $\tau$.

***Definition 9:*** For a collection of lattices $\tau$ having a unit co-volume in Euclidean space $R^n$ i.e. $Vol\left(\frac{R^n}{L}\right) = 1$. We define the Hermite Constant $\gamma_n$ for any positive integer $n$ as the square of the maximum of $\lambda_1(\tau)$ which is the least length of a nonzero element of $\tau$ as defined earlier over all the lattices in that collection. The Hermite constant can be alternatively defined for dimension $n$ by Lionnais[20] as follows:

$$\gamma_n = \frac{sup_f min_{x_i} f(x_i, x_i, \ldots x_n)}{[discriminant(f)]^{1/n}}$$

It can also be alternatively defined as

$$\gamma_n = (\delta_n/v_n)^{2/n}$$

Here $v_n$ is the volume of the $n$-dimensional sphere.

Finding the exact value of $\gamma_n$ is a very difficult problem; Martinet[21] gave the exact values for $1 \leq n \leq 8$. The values of $\gamma_n{}^n$ are $1, 4/3, 2, 4, 8, 64/3, 64, 256, \ldots$ (OEIS[22] A007361 and A007362, Gruber and Lekkerkerker[23], Weisstein[24]). But there are various approximation and very tight bounds that are found by mathematicians over the years, some of which we list here. But first we derive a simple one which follows immediately from the Minkowski's Convex Body Theorem as proved earlier in **Theorem 3.**

Note that **Corollary 4** of **Theorem 3** can be restated in other words as

$$\gamma_d \leq \left(\frac{4}{v_d}\right)^{\frac{2}{d}}, d \geq 1$$

The formula of a $n-$dimensional volume of a Euclidean ball of radius $R$ is[25]:

$$v_d(R) = \frac{\pi^{d/2}}{\vartheta(\frac{d}{2}+1)} R^n$$

From this two previous equations Hans Frederick Blichfeldt[26,27] gave a stronger estimate of $\gamma_d$ as follows:

$$\gamma_n \leq \left(\frac{2}{\pi}\right)\vartheta\left(2 + \frac{n}{2}\right)^{\frac{2}{n}}$$

Here $\vartheta(x)$ is the gamma function. Kitaoka[27] gave an upper bound to the Hermite constant in $n$ dimensions by stating the inequality $\gamma_n \leq \left(\frac{4}{3}\right)^{\frac{n-1}{2}}$. Now we also know that $v_d$ which denotes the volume of the closed unitary hyperball of $R^d$ can be approximated as

$$v_d(R) = \frac{\pi^{d/2}}{\vartheta(\frac{d}{2}+1)} R^n \sim \left(\frac{2e\pi}{n}\right)^{n/2} \frac{1}{\sqrt{\pi n}}$$

From this approximation we can approximate $\lambda_1$ to be $(vol(\tau)/v_d)^{1/d} \approx \sqrt{d/2\pi e}\, vol(\tau)^{1/d}$. Therefore we have $\gamma_n = \lambda_1{}^2/vol(\tau)^{2/d} \approx d/2\pi e$. Mathematicians have proven rather very tight asymptotical bounds for Hermite's constant. One lower and upper bound states that

$$\frac{n}{2\pi e} + \frac{\log(\pi n)}{2\pi e} + o(1) \leq \gamma_n \leq \frac{1.744 d}{2\pi e}(1 + o(1))$$

From this we can derive an inequality for large values of $n$ as derived by Weisstein[24], Weisstein[24] has given the upper and lower bounds of the Hermite Constant $\gamma_n$ as:

$$\frac{1}{2\pi e} \leq \frac{\gamma_n}{n} \leq \frac{1.744..}{2\pi e}$$

We shall now discuss about the Minkowski-Hlawaka which is stated without proof by Minkowski[28] and proved later by Hlawka[29].

***Theorem 7[Lattice packing and the Minkowski-Hlawka theorem]:*** It states that there is a lattice of dimension $n$ in Euclidean space, such that the best packaging of hyperspheres with centers at the lattice points has density $\rho$ satisfying the following inequality

$$\rho \geq \zeta(n)/2^{n-1}$$

Here $\zeta(n)$ is the Riemann zeta function. The proof of this theorem, however, is non-constructive and how to really construct packings that are this dense is still not understood[30].

The following generalization of the Minkowski–Hlawka theorem was proven by Siegel[31].

***Theorem 8:*** The average number of nonzero lattice vectors in S which is a bounded set in the Euclidean space of dimension $n$ with Jordan volume $vol(S)$ is $vol(S)/D$. Here all lattices with a volume $D$ fundamental domain are taken over by the average, Similarly, $vol(S)/D\zeta(n)$ is the average number of primitive lattice vectors in $S$.

## 5. VORONOÏ CELL, PACKING RADIUS AND COVERING RADIUS

In the previous section we have studied about successive minima of any given full-rank lattice $\tau$ of a given rank $n$ and derived some very interesting properties of them. In this section we shall discuss yet another other useful parameters of lattice. Those are the Voronoï Cell, packing radius and the covering radius of the lattice. We shall start with the definition Voronoï Cell followed by that of packing and covering radius. This section about Voronoï Cell is composed in the form of the paper originally written by Moustrou[17].

***Definition 8:*** The Voronoï cell of a lattice $\tau \subset R^n$ is a region associated to a given full-rank lattice $\tau$ which consists of the points of $R^n$ that is closer to 0 than to any other vector of $\tau$

$$V = V_A = \{z \in R^n | such\ that\ x \in \tau, \|z - x\| \geq \|z\| \}$$

and it is a fundamental region of $\tau$. The Voronoï Cell is defined by the Voronoï Vectors of the given lattice $\tau$. A vector $v$ belonging to the lattice $\tau$ is a Voronoï Vector of the same lattice if the intersection of Voronoï Cell $V$ and an n-dimensional hyper plane defined as

$$H_v = \{x \in R^n | \langle x, v \rangle = \frac{1}{2} \langle v, v \rangle \}$$

Contains some elements i.e. it is not empty. The vector $v$ is said to be a relevant vector if this intersection $H_v$ defined earlier is a $(n-1)$ dimensional hyper plane of the Voronoi Cell. Note that the relevant vectors in itself enough to fully describe the Voronoi Cell $V$, a vector that belongs to Voronoi Cell $V$ if and only if for all relevant Voronoï vector $v$ we have the following inequality holding true

$$|\langle x, v \rangle| \leq \frac{1}{2} \langle v, v \rangle$$

Another proposition composed originally by Moustrou[17] states that a non-trivial vector $v$ is a Voronoï vector of the full-rank lattice $\tau$ if it is the shortest vector in the coset $v + 2\tau$. Moreover, $v$ is relevant if and only if $\pm v$ are the only shortest vectors in that coset. A nice treatise on the characterization and properties of the Voronoï vectors and Cell can be found by Conway and Sloane.[18]

This section is composed in the form of the paper originally written by Micciancio[6]

***Definition 9:*** The packing radius $r(\tau)$ of any given full-rank lattice $\tau$ of a given rank $n$ is defined as the largest radius $r > 0$ such that the open balls $B(v,r) = \{x: \|x - v\| < r\}$ are centered on all lattice points do not intersect.

***Definition 10:*** The covering radius $\mu(\tau)$ of any given full-rank lattice $\tau$ of a given rank $n$ is defined as the smallest radius $r > 0$ such that the closed balls $\bar{B}(v,r) = \{x: \|x - v\| \leq r\}$ shall cover the entire space or $span(\tau)$ is a subset of $\bigcup_{v \in \tau} \bar{B}(v, \mu)$

We shall now state and prove some relationships between the packing and covering radius in terms of bounds in the next few theorems. The first simple and elegant proof is due to Oded Regev[20].

***Theorem 5:*** For a given lattice $\tau$, the covering radius $\mu(\tau)$ is at least greater than $1/2\lambda_n$

***Proof:*** Let us consider the open ball $B(0, \lambda_n)$. Now note that all the lattice points inside this open ball should be in a $(n − 1)$ dimensional plane since every lattice point having length lesser than $\lambda_n$ must be in an $(n − 1)$ dimensional plane. Now perpendicular to the hyper plane we go a distance $1/2\lambda_n$ to get a point $x$. By the construction itself any lattice point outside the ball and also inside the ball( which lies in the (n-1) dimensional hyper plane ) must be at least $1/2\lambda_n$ distance from $x$ i.e the covering radius is at least greater than $1/2\lambda_n$.

***Theorem 6:*** Show that for any $n -$dimensional lattice $\tau$, the covering radius is at most $\frac{\sqrt{n}}{2}\lambda_n$

***Proof:*** Let us construct a sub-lattice from $n$ linearly independent lattice vectors $\{a_i\}$. of length at most $\lambda_n$. We shall bound using the fundamental region $C(A^*)$ the covering radius of the sub-lattice $L(A)$ generated from lattice vectors $\{a_i\}$. Therefore the cube of side length $\mu(\tau)$ is contained inside the ball of the radius $\frac{\sqrt{n}}{2}\lambda_n$, as we did earlier in Corollary 1 of Minkowski's convex body theorem. Thus we can write $\mu(\tau) \leq \frac{\sqrt{n}}{2}\lambda_n$

The next theorem gives us another lower limit to the covering radius of the lattice $\tau$, before proving the theorem note that we can use Minkowski's convex body theorem to get an upper bound on the packing radius of the lattice $\frac{1}{2}\lambda_n \leq (det(\tau)/V_n)^{\frac{1}{n}}$. (Here $V_n$ is the volume of the ball of radius 1 in $R^n$). We can use a similar proof using an argument similar to the proof of Blichfeldt theorem to prove a lower bound on the covering radius. We can prove that the covering radius is at least $\mu(\tau) \geq (det(\tau)/V_n)^{\frac{1}{n}}$, where $V_n$ is the volume of the unit ball in n-dimensional ***R*** space.

## 6. APPLICATION

Lattices play a very fundamental role in the theory of numbers.

Now we present some examples of number theoretic application of Minkowski's theorem. We first start with Fermat theorem on sum of two squares which states that all primes which can be written as 4k + 1 can be written as a sum of two squares. The result was proved by Axel Thue using the pigeonhole principle. We now present a proof using Minkowski's lattice point theorem originally derived by Jiri Matousek[7] in his book Lectures on Discrete Geometry and Dong[21].

***Lemma 4:*** If p is a prime with $p \equiv 1 \ (mod \ 4)$ then -1 is a quadratic residue modulo p.

***Proof:*** Note that the equation $i^2 = 1$ has two solutions, namely $i = 1 \ and \ i = -1$ in the field $F$. Hence for any $i \neq \pm 1$ there will exists exactly one element such that if $j \neq i$, $ij = 1$ and therefore we can divided into pairs all the elements present in $F^* \setminus \{-1, 1\}$ in which all the product of elements coming in pair will be equal to 1. Therefore, $(p − l)! = 1 \cdot 2 \cdots (p − l) \equiv -1 \ (mod p)$. Now by contradiction, let us assume that the equation $i^2 = $ -1 has no solution in F. Then all the elements of F* can be divided into pairs such that the product of the elements in each pair will be -1. But we have $(p − 1)/2$ pairs, which is an even number. The product of which will be 1 as shown below.

Hence $(p − l)! \equiv (-1)^{(p-1)/2} = 1$ ; we therefore arrive at a contradiction.

***Theorem 7:*** Primes of the form 4k + 1 can be expressed as a sum of two squares

***Proof:*** Let us assume $p$ to be a prime of the form $4k + 1$. Now as proven earlier in Lemma 4 we know that if p is a prime with $p \equiv 1 \ (mod \ 4)$ then -1 is a quadratic residue modulo p.

Now let us assume $q^2 \equiv -1 (mod \ p)$. Let us also take two vectors $z_1 = (1, q)$ and $z_2 = (0, p)$. We define the lattice $\tau$ as the span of the basis vectors $\{z_1, z_2\}$. The volume of the lattice shall be $det(\tau) = p$. Now consider a curve defined as:

# Lattices and the Geometry of Numbers

$$C = \{(x,y): x^2 + y^2 < 2p\} \text{ in } R^2$$

Now according to Minkowski Theorem we discussed earlier if $Vol(X) > 2^n Vol(\tau)$ then $X$ must contain a nonzero lattice point $\gamma = (a,b) \in \tau$.

$$Vol(C) = 2\pi p > 4p = 2^2 \det(\tau)$$

Now note that $(a,b) = iz_1 + jz_2 = (i, iq + jp) \in R^2$, We can therefore say that

$$a^2 + b^2 = i^2 + (iq + jp)^2 \equiv (q^2 + 1)i^2 \equiv 0 (mod\ p)$$

But as we have defined earlier $a^2 + b^2$ must be lesser than $2p$, the only way the above condition holds is when $a^2 + b^2 = p$. Hence we proved that primes of the form $4k + 1$ can be expressed as a sum of two squares.

Let's see another beautiful application of Minkowski's Theorem by proving Dirichlet's theorem on Diophantine approximation. Let us assume $\alpha$ to be any given real number. We now wish to approximate $\alpha$ with rational numbers which are numbers that can be expressed in the form of $p/q$. Now notice that we can approximate it with any given precision if we continue to expand the continued fraction of it, but to make it bit more challenging we also put a constraint denominator of a rational number; it should not exceed a given integer $Q$. We now state and prove the theorem as done by Matousek[7].

***Theorem 8(Dirichlet):*** Let $\alpha$ be a real number and $Q$ a positive integer. Then there are integers $p$ and $q$ with $0 < q \leq Q$ and $|\alpha - p/Q| \leq 1/Q$.

***Proof:*** Let us consider the volume $C$ in $R^2$ that is closed by the lines

$$y \leq \alpha x + \frac{1}{Q}, y \geq \alpha x - \frac{1}{Q}, x \leq Q, x \geq -Q.$$

The volume of the parallelopoid will be $vol(C) = 4Q * \left(\frac{1}{Q}\right) = 4$

The following theorem was first proved by Lagrange and states that every positive integer can be written as the sum of four squares of integers. It can also be derived with the help of Minkowski's Theorem. We shall state and prove the theorem as done by Shmonin[22].

***Theorem 9(Lagrange):*** For every positive integer x, there are integers $x_1, x_2, x_3, x_4$ such that $x = x_1^2 + x_2^2 + x_3^2 + x_4^2$.

***Proof:*** This theorem needs to be proved only for integers other than the prime number because if it is not a prime number, then we can factorize the number to at least two prime factors let's say $x = x_1^2 + x_2^2 + x_3^2 + x_4^2$ and $y = y_1^2 + y_2^2 + y_3^2 + y_4^2$, then the product $xy$ can be written in the following form:

$xy = (x_1y_1 + x_2y_2 + x_3y_3 + x_4y_4)^2 + (x_1y_2 - x_2y_1 + x_3y_4 - x_4y_3)^2 + (x_1y_3 - x_3y_1 + y_4 + x_4y_3)^2 + (x_1y_4 + x_2y_3 - x_3y_2 - x_4y_1)^2$.

Before proving the theorem let us show that if x is prime integer there are must exist two integers $y$ and $z$ such that $y^2 + z^2 + 1 \equiv 0\ (mod\ x)$, now let us consider 2 sets of numbers defined as follows:

$$S_1 := \{y^2 mod x : 0 \leq y \leq (x-1)/2\} \text{ And } S_2 := \{-z^2 - 1\ mod x : 0 \leq z \leq (x-1)/2\}$$

It is easy to see that for $x = 2$ this becomes a trivial case, otherwise it should be a odd integer. Also notice that for different $0 \leq y_1, y_2 \leq (x-1)/2$ the inequality $(y_1 mod\ x) \neq (y_2 mod\ x)$ holds. Therefore the cardinality of $S_1$ is $(x + 1)/2$. Similarly it can be shown that the cardinality of $S_2$ is $(x + 1)/2$. Therefore we can say that there must exists an intersection such that $y^2 \equiv -(z^2 + 1) \pmod{x}$, i.e $y^2 + z^2 + 1 \equiv 0 \pmod{x}$.

We shall now define a lattice $\tau$ whose basis vector matrix $B$ is as written below:

$$B = \begin{bmatrix} x & 0 & y & z \\ 0 & x & z & -y \\ 0 & 0 & 1 & 0 \\ 0 & 0 & 0 & 1 \end{bmatrix}$$

Let us now apply Minkowski's theorem to the lattice $\tau$. By Minkwoski's convex body theorem we know that if $Vol(X) > 2^n det(\tau)$ then the $X$ must have a non-trivial point inside it.

Now note that $det(\tau) = x^2$. For the application of Minkwoski's convex body theorem let us also consider the ball which is defined as

$$B = \{[x_1, x_2, x_3, x_4]: x_1^2 + x_2^2 + x_3^2 + x_4^2 < 2x\}$$

The Volume of the ball shall be

$$vol(B) = \left(\frac{1}{2}\right) \pi^2 \sqrt{2x}^4 = 2\pi^2 x^2 > 2^4 x^2$$

Therefore there must exist a non-trivial point $[x_1, x_2, x_3, x_4]$, inside the ball $B$ Therefore we can write

$$x_1^2 + x_2^2 + x_3^2 + x_4^2 < 2x$$

But we can also derived that

$$x_1^2 + x_2^2 + x_3^2 + x_4^2 = (x\lambda_1 + y\lambda_3 + z\lambda_4)^2 + (x\lambda_2 + z\lambda_3 - y\lambda_4)^2 + \lambda_3^2 + \lambda_4^2 \equiv (1 + y^2 + z^2)(\lambda_3^2 + \lambda_4^2) \pmod{x} \equiv 0 \pmod{x}.$$

Since every vector $[x_1, x_2, x_3, x_4] = B[\lambda_1, \lambda_2, \lambda_3, \lambda_4]^T$ for the given lattice $\tau$. The above two conditions will hold simultaneously only when $x = x_1^2 + x_2^2 + x_3^2 + x_4^2$ which proves our theorem for prime numbers also.

## 7. CONCLUSION

In this paper we have reviewed all the major developments that have occurred in the field of geometry of numbers. In this paper we have also given a broad overview of the concept of lattice and then discuss about the geometrical properties it has starting from the Minkowski convex body theorem which led to the beginning of the new field of geometry of numbers. We have also saw applications to the theory of numbers through suitable examples.